\documentclass{article}
\usepackage{amsmath,amssymb,amsthm}

\newtheorem{theorem}{Theorem}
\newtheorem{proposition}{Proposition}
\begin{document}
\title{Universal Metric Spaces According to W.~Holsztynski}
	\author{Sergei~Ovchinnikov \\
	Mathematics Department\\
	San Francisco State University\\
	San Francisco, CA 94132\\
	sergei@sfsu.edu} 
\date{\today}
\maketitle

\begin{abstract}
\noindent
In this note we show, following W.~Holsztynski, that there is a continuous metric $d$ on $\mathbb{R}$ such that any finite metric space is isometrically embeddable into $(\mathbb{R},d)$.
\end{abstract}

\bigskip
Let $\mathcal{M}$ be a family of metric spaces. A metric space $U$ is said to be a \emph{universal space} for $\mathcal{M}$ if any space from $\mathcal{M}$ is (isometrically) embeddable in $U$.

Fr\'{e}chet~\cite{mF10} proved that $\ell^\infty$ (the space of all bounded sequences of real numbers endowed with the $\sup$ norm) is a universal space for the family $\mathcal{M}$ of all separable metric spaces. Later, Uryson~\cite{pU27} constructed an example of a separable universal space for this $\mathcal{M}$ ($\ell^\infty$ is not separable). 

In this note we establish the following result which is Theorem~5 in~\cite{wH78}.

\begin{theorem}
There exists a metric $d$ in $\mathbb{R}$, inducing the usual topology, such that every finite metric space embeds in $(\mathbb{R},d)$.
\end{theorem}

Our proof essentially follows the original Holsztynski's approach~\cite{wH00}.

We say that a metric space $(X,d)$ is $\varepsilon$--\emph{dispersed} if $d(x,y)\geq\varepsilon$ for all $x\not=y$ in $X$ $(\varepsilon>0)$. Clearly, any $\varepsilon$-- metric space is also $\varepsilon'$--dispersed for any positive $\varepsilon'<\varepsilon$. The following proposition will be used to construct universal spaces for particular families $\mathcal{M}$.

\begin{proposition} \label{Proposition1}
Let $f:X\rightarrow Y$ be a continuous surjection from a metric space $(X,d)$ onto a metric space $(Y,D)$. Then $(X,d_\varepsilon)$ where
\begin{equation} \label{Equation1}
d_\varepsilon(x,y) = \max\{\min\{d(x,y),\varepsilon\},D(f(x),f(y))\}
\end{equation}
is a universal space for the family of $\varepsilon$--dispersed subspaces of $(Y,D)$ and metrics $d$ and $d_\varepsilon$ are equivalent on $X$.
\end{proposition}

\begin{proof}
$d_\varepsilon$ is a distance function. Indeed, $d_\varepsilon$ is symmetric and $d_\varepsilon(x,y)=0$ if and only if $x=y$. We have
\begin{align*}
\max\{&\min\{d(x,y),\varepsilon\},D(f(x),f(y))\}+\max\{\min\{d(y,z),\varepsilon\},D(f(y),f(z))\}\geq \\
&\geq D(f(x),f(y))+D(f(y),f(z))\geq D(f(x),f(z))
\end{align*}
and
\begin{align*}
\max\{&\min\{d(x,y),\varepsilon\},D(f(x),f(y))\}+\max\{\min\{d(y,z),\varepsilon\},D(f(y),f(z))\}\geq \\
&\geq \min\{d(x,y),\varepsilon\}+\min\{d(y,z),\varepsilon\}=\min\{d(x,y)+d(y,z),\varepsilon\}\geq \\
&\geq \min\{d(x,z),\varepsilon\}
\end{align*}
Hence, $d_\varepsilon(x,y)+d_\varepsilon(y,z)\geq d_\varepsilon(x,z)$.

Let $Z$ be an $\varepsilon$--dispersed subspace of $Y$. Since $f$ is surjective, for any $z\in Z$, there is $x_z\in X$ such that $f(x_z)=z$. Let $X'=\{x_z:z\in Z\}$. By (\ref{Equation1}), $d_\varepsilon(x,y)=D(f(x),f(y))$ for all $x,y\in X'$. Thus $f$ establishes an isometry between $(Z,D)$ and $(X',d_\varepsilon)$.

\end{proof}

In what follows, $\mathcal{M}$ is the family of all finite metric spaces.

We define
\begin{equation*}
I^n = \{\bar{x}=(x_1,\ldots,x_n)\in\mathbb{R}^n:0\leq x_i\leq n,\;1\leq i\leq n \}
\end{equation*}
and $J_n=[n-1,n]$ for $n\geq 1$. $I^n$ is a metric space with the distance function
\begin{equation*}
D(\bar{x},\bar{y}) =\max\{|x_i-y_i|:1\leq i\leq n\},
\end{equation*}
and $J_n$ is a metric space with the usual distance.

Let $(X,d)$ be a finite metric space. We define 
$$
p=|X|,\;q=\lfloor\text{Diam}(X)\rfloor,\;r=\lfloor\varepsilon^{-1}\rfloor,
$$
where $\varepsilon=\min\{d(x,y):x\not= y\}$, and $n=\max\{p,q,r\}$. Clearly, $(X,d)$ is $\frac{1}{n}$--dispersed.

\begin{proposition} \label{Proposition2}
\emph{(The Kuratowski embedding)} $(X,d)$ is embeddable into $I^n$.
\end{proposition}

\begin{proof}
Let $X=\{x_1,\ldots,x_p\}$. We define $f:X\rightarrow I^n$ by
\begin{equation*}
f(x_i) = (d(x_i,x_1),\ldots,d(x_i,x_p),\underbrace{0,\ldots,0}_{n-p\;0\text{'s}}),
\end{equation*}
for $1\leq i\leq p$. (Since $n\geq \text{Diam}(X)$,\;$f(x_i)\in I^n$.) We have, by the triangle inequality,
\begin{equation*}
D(f(x_k),f(x_m))=\max_j\{|d(x_k,x_j)-d(x_m,x_j)|\}\leq d(x_k,x_m).
\end{equation*}
On the other hand, $|d(x_k,x_j)-d(x_m,x_j)|=d(x_k,x_m)$ for $j=m$. Therefore, $D(f(x_k),f(x_m))=d(x_k,x_m)$ for all $1\leq k,m\leq p$.

\end{proof}

Let $f_n$ be a continuous surjection from $J_n$ onto $I^n$ (a ``Peano curve''~\cite[IV(4)]{jD89}) and $d_n(x,y)$ be the distance function on $J_n$ defined by~(\ref{Equation1}) for $\varepsilon =\frac{1}{n}$. Note, that $d_n$ is equivalent to the usual distance on $J_n$. By Proposition~\ref{Proposition1}, $J_n$ is a universal space for any $\frac{1}{n}$--dispersed subspace of $I^n$. By Proposition~\ref{Proposition2}, any finite metric space is embeddable in $(J_n,d_n)$ for some $n$.

It is easy to show that there is a continuous distance function on $\mathbb{R}$ that coincides with $d_n(x,y)$ on $J_n$ for all $n$. Indeed, let $d_1(x,y)$ and $d_2(x,y)$ be two continuous distance functions on intervals $[a,b]$ and $[b,c]$, respectively. Then $d(x,y)$ defined by
\begin{equation*}
d(x,y) = \begin{cases}
	d_1(x,y),	&\text{if $x,y\in [a,b]$,} \\
	d_2(x,y),	&\text{if $x,y\in [b,c]$,} \\
	d_1(x,b)+d_2(b,y),	&\text{if $x\in [a,b]$ and $y\in [c,d]$.}
\end{cases}
\end{equation*}
is a continuous distance function on $[a,c]$. In fact, thus defined $d$ is equivalent to the usual metric on $[a,c]$. 

By applying this process consecutively, we obtain a required distance function on $\mathbb{R}$.

\end{document}